\documentclass[12pt,twoside]{article}

\begin{document}
\centerline{}

\centerline {\Large{\bf Fuzzy ideals in $\Gamma-$semiring}}

\centerline{}

\newcommand{\mvec}[1]{\mbox{\bfseries\itshape #1}}

\centerline{\bf {Jayanta Ghosh, T.K. Samanta }}

\centerline{}

\centerline{Department of Mathematics,}
\centerline{Manickpur Adarsha Vidyapith, West Bengal, India. }
\centerline{e-mail: ghoshjay$_{-}$04@yahoo.com}

\centerline{Department of Mathematics, Uluberia
College, West Bengal, India.}
\centerline{e-mail: mumpu$_{-}$tapas5@yahoo.co.in}
\centerline{}

\newtheorem{Theorem}{\quad Theorem}[section]

\newtheorem{definition}[Theorem]{\quad Definition}
\newtheorem{proposition}[Theorem]{\quad proposition}

\newtheorem{theorem}[Theorem]{\quad Theorem}

\newtheorem{remark}[Theorem]{\quad Remark}

\newtheorem{corollary}[Theorem]{\quad Corollary}

\newtheorem{note}[Theorem]{\quad Note}
\newtheorem{lemma}[Theorem]{\quad Lemma}

\newtheorem{example}[Theorem]{\quad Example}

\newtheorem{notation}[Theorem]{\quad Notation}
\begin{abstract}
In this paper we have studied the relation between the fuzzy left (respectively right) ideals of $\Gamma-$semiring and that of operator semiring.
Thereafter, we have established that the Lattices of all fuzzy left (respectively right) ideal of $\Gamma-$semiring is equivalent to that of
Left operator semiring. Also we have established few properties relating the $k-$ideals and $h-$ideals of $\Gamma-$semiring with that of
operator semiring.
\end{abstract}
\\\\
\section{Introduction}
The notion of $\Gamma$ in algebra was first introduced by N. Naobuswa\cite{Nobuswa}
in 1964 and also he defined the $\Gamma-$ring.
In 1969, J. Luh\cite{Luh} introduced the concept of left operator ring and right
operator ring of $\Gamma-$ring.
In 1995, M. M. K. Rao\cite{MMK} introduced the concept of
$\Gamma-$semiring as a generalization of semiring and $\Gamma-$ring. Thereafter
S. K. Sardar and T. K. Dutta\cite{Sardar} modified the definition of
$\Gamma -$semiring of Rao\cite{MMK} and then they redefined the left operator
semiring and right operator semiring of a $\Gamma-$semiring and obtained a few
interesting properties. Later on, much has been developed on this concepts by different
researchers.\\
Fuzzy set theory was first introduced by Zadeh\cite{zadeh} in 1965
and thereafter several authors contributed different articles on
this concept and applied it on different branches of pure and applied
mathematics. In 1971, Rosenfeld\cite{Rosenfeld} defined fuzzy subgroups, fuzzy ideals and studied its important properties.
Thereafter in 1982, Liu\cite{Liu} introduced and developed basic results of fuzzy subrings and fuzzy ideals of a ring.
In 1992, Jun and Lee\cite{Jun} introduced the notion of fuzzy ideal in $\Gamma-$ring and studied a few properties. In 2005,
Dutta and Chanda\cite{Dutta} studied the structures of fuzzy ideals of $\Gamma-$ring via operator rings of $\Gamma-$ring.\\
In this paper, we have established a few results in respect of fuzzy left(respectively right) ideals, fuzzy ideals of a $\Gamma-$semiring and
its operator semirings. If $\sigma$ is a fuzzy ideal of a $\Gamma-$semiring then we have proved $\sigma^{+}\,^{'}$ is a fuzzy ideal of the corresponding left operator semiring. Also, if $\mu$ is a fuzzy ideal of left operator semiring then $\mu^{+}$ is a fuzzy ideal of the corresponding $\Gamma-$semiring.
Then it is shown that there exist an inclusion preserving bijection $\sigma\,\rightarrow\,\sigma^{+}\,^{'}$ between
the Lattices of all fuzzy right ideals (respectively fuzzy ideals)
of a $\Gamma-$semiring and the Lattices of all fuzzy right ideals (respectively fuzzy ideals) of the corresponding left operator semiring.
Similarly the above results hold for right operator semiring of a $\Gamma-$semiring.
Also we have studied similar results relative to fuzzy $k-$ideals, fuzzy $h-$ideals
of a $\Gamma-$semiring and its operator semirings.

\section{Preliminaries}
This section contain some basic definitions and preliminary results which will be needed in the sequel.\\
\begin{definition}\cite{Sardar}
Let $S$ and $\Gamma$ be two additive commutative semigroups. Then $S$ is called a $\Gamma$$-$semiring if
there exists a mapping $S \times\Gamma\times S\rightarrow S$ (image to be denoted by $a\alpha b$ where $a,b\in S$ and $\alpha\in\Gamma$) satisfying the following conditions:\\
(1) $a\,\alpha\,(b+c)\;=\;a\,\alpha\,b\;+\;a\,\alpha\,c$\\
(2) $(a+b)\,\alpha\,c\;=\;a\,\alpha\,c\;+\;b\,\alpha\,c$\\
(3) $a\,(\alpha+\beta)\,c\;=\;a\,\alpha\,c\;+\;a\,\beta\,c$\\
(4) $a\,\alpha\,(b\,\beta\,c)\;=\;(a\,\alpha\,b)\,\beta\,c$\\
for all $a,b,c\,\in\,S$ and for all $\alpha,\,\beta\,\in\,\Gamma$.
\end{definition}
\begin{definition}
Let $S$ be a $\Gamma-$semiring and  $\mu$ be a fuzzy subset of $S$.
Then $\mu$ is called a fuzzy left ideal of $S$ if\\
(1) $\mu(a+b)\,\geq\,min(\,\mu(a)\,,\,\mu(b)\,)$\\
(2) $\mu(a\,\alpha\,b)\,\geq\,\mu(b)$\\
for all $a,b\,\in\,S$ and for all $\alpha\,\in\,\Gamma$.
\end{definition}
\begin{definition}
Let $S$ be a $\Gamma-$semiring and $\mu$ be a fuzzy subset of $S$.\\
Then $\mu$ is called a fuzzy right ideal of $S$ if\\
(1) $\mu(a+b)\,\geq\,min(\,\mu(a)\,,\,\mu(b)\,)$\\
(2) $\mu(a\,\alpha\,b)\,\geq\,\mu(a)$\\
for all $a,b\,\in\,S$ and for all $\alpha\,\in\,\Gamma$.
\end{definition}
Note: If a fuzzy subset $\mu$ of $\Gamma-$semiring $S$ is both fuzzy left and fuzzy right ideal of $S$ then $\mu$ is called a fuzzy two-sided ideal or simply fuzzy ideal of $S$.\\
\begin{definition}
Let $S$ be a $\Gamma-$semiring. A fuzzy ideal $\mu$ of $S$ is called fuzzy $k-$ideal of $S$ if $\;\mu(x)\,\geq\,min(\,\mu(x+y)\,,\,\mu(y)\,)\hspace{1cm}$for all $x,y\,\in\,S$.
\end{definition}
\begin{definition}
Let $S$ be a $\Gamma$$-$semiring. A fuzzy ideal $\mu$ of $S$ is called fuzzy $h-$ideal of $S$ if for all $x,z,y_{1},y_{2}\,\in\,S$ s.t. $x+y_{1}+z=y_{2}+z$ implies $\mu(x)\,\geq\,min(\,\mu(y_{1})\,,\,\mu(y_{2})\,)$.
\end{definition}
Similarly we define fuzzy one-sided $k-$ideal and one-sided $h-$ideal of $S$.
\begin{definition}\cite{Sardar}
\textbf{Left operator semiring and Right operator\\ semiring of a $\Gamma$-semiring}
\\Let $S$ be a $\Gamma$-semiring and $F$ be the free addtive commutative semigroup generated by $S\times\Gamma$. Then the relation $\rho$ on $F$  defined by
$\sum_{i}\,(x_{i},\alpha_{i})\,\rho\,\sum_{j}\,(y_{j},\beta_{j})$ iff $\sum_{i}\,x_{i}\alpha_{i}a\,=\,\sum_{j}\,y_{j}\beta_{j}a$ for all $a\,\in\,S$,
is congruence on $F$. We denote the congruence class containing $\sum_{i}\;(x_{i},\alpha_{i})$ by $\sum_{i}\;[x_{i},\alpha_{i}]$. Then $F/\rho$  is an
additive commutative semigroup. Now we define a multiplication on $F/\rho$  by $(\sum_{i}\;[x_{i},\alpha_{i}])(\sum_{j}\,[y_{j},\beta_{j}])\,=\,\sum_{i},_{j}\,[x_{i}\alpha_{i}y_{j},\beta_{j}]$. Then $F/\rho$ forms a semiring with multiplication defined above. We denote this semiring by $L$ and call it the left operator semiring of the $\Gamma$-semiring $S$.\\
Dually we define the right operator semiring $R$ of the $\Gamma$-semiring $S$ where\\$R\,=\,\{\,\sum_{i}\,[\alpha_{i},x_{i}]\,:\,\alpha_{i}\,\in\,\Gamma
,\,x_{i}\,\in\,S\,\}$ and the multiplication on $R$ is defined as $(\sum_{i}\,[\alpha_{i},x_{i}])(\sum_{j}\,[\beta_{j},y_{j}])\,=\,\sum_{i},_{j}\,[\alpha_{i},x_{i}\beta_{j}y_{j}]$.\\\\
We also note here that for $[x+y,\alpha]\,,\,[x,\alpha+\beta]\,\in\,L$,\\$[x+y,\alpha]\,=\,[x,\alpha]+[y,\alpha]$ and $[x,\alpha+\beta]\,=\,[x,\alpha]+[x,\beta]$.\\Similarly for $[\alpha,x+y]\,,\,[\alpha+\beta,x]\,\in\,R$,\\
$[\alpha,x+y]\,=\,[\alpha,x]+[\alpha,y]$ and $[\alpha+\beta,x]\,=\,[\alpha,x]+[\beta,x]$.
\end{definition}
\begin{definition}
Let $S$ be a $\Gamma-$semiring and $L$ be its left operator semiring. Let $\mu$ be a fuzzy subset of $L$. Then $\mu$ is called a fuzzy left ideal of $L$ if\\
(1) $\mu(\;\sum_{i}\,[x_{i},\alpha_{i}]\,+\,\sum_{j}\,[y_{j},\beta_{j}]\;)\;
\geq\;min(\;\mu(\,\sum_{i}\,[x_{i},\alpha_{i}]\,)\,,\,\mu(\,\sum_{j}\,[y_{j},\beta_{j}]\,)\;)$\\
(2) $\mu(\;\sum_{i}\,[x_{i},\alpha_{i}]\,.\,\sum_{j}\,[y_{j},\beta_{j}]\;)\;\geq\;\mu(\,\sum_{j}\,[y_{j},\beta_{j}]\,)$\\
for all $\sum_{i}\,[x_{i},\alpha_{i}]\,,\,\sum_{j}\,[y_{j},\beta_{j}]\,\in\,L$.
\end{definition}
\begin{definition}
Let $S$ be a $\Gamma-$semiring and $L$ be its left operator semiring. Let $\mu$ be a fuzzy subset of $L$. Then $\mu$ is called a fuzzy right ideal of $L$ if\\
(1) $\mu(\;\sum_{i}\,[x_{i},\alpha_{i}]\,+\,\sum_{j}\,[y_{j},\beta_{j}]\;)\;
\geq\;min(\;\mu(\,\sum_{i}\,[x_{i},\alpha_{i}]\,)\,,\,\mu(\,\sum_{j}\,[y_{j},\beta_{j}]\,)\;)$\\
(2) $\mu(\;\sum_{i}\,[x_{i},\alpha_{i}]\,.\,\sum_{j}\,[y_{j},\beta_{j}]\;)\;\geq\;\mu(\,\sum_{i}\,[x_{i},\alpha_{i}]\,)$\\
for all $\sum_{i}\,[x_{i},\alpha_{i}]\,,\,\sum_{j}\,[y_{j},\beta_{j}]\,\in\,L$.
\end{definition}
Note: If a fuzzy subset $\mu$ of $L$ is both fuzzy left and fuzzy right ideal of $L$ then $\mu$ is called fuzzy two-sided ideal or simply fuzzy ideal of $L$.
\begin{definition}
Let $S$ be a $\Gamma$-semiring and L be its left operator semiring. A fuzzy ideal $\mu$ of L is called fuzzy $k$-ideal if for all $\sum_{i}\,[x_{i},\alpha_{i}]\,,\,\sum_{j}\,[y_{j},\beta_{j}]\,\in\,L$,\\
$\mu(\,\sum_{i}\,[x_{i},\alpha_{i}]\,)\,\geq\,min(\,\mu(\,\sum_{i}\,[x_{i},\alpha_{i}]+\sum_{j}\,[y_{j},\beta_{j}]\,)\,,\,\mu(\sum_{j}\,[y_{j},\beta_{j}])\,)$.
\end{definition}
\begin{definition}
Let $S$ be a $\Gamma$-semiring and $L$ be its left operator semiring. A fuzzy ideal $\mu$ of $L$ is called fuzzy $h$-ideal if for all
$\sum_{i}\,[x_{i},\alpha_{i}]\,,\,\sum_{j}\,[y_{j},\beta_{j}],\\\,\sum_{k}\,[z_{k},\gamma_{k}],\,\sum_{l}\,[u_{l},\delta_{l}]\,\in\,L$
s.t. $\sum_{i}\,[x_{i},\alpha_{i}]+\sum_{j}\,[y_{j},\beta_{j}]+\sum_{k}\,[z_{k},\gamma_{k}]\,=\,\sum_{l}\,[u_{l},\delta_{l}]+\sum_{k}\,[z_{k},\gamma_{k}]$ implies\\$\mu(\,\sum_{i}\,[x_{i},\alpha_{i}]\,)\,\geq\,min(\,\mu(\sum_{j}\,[y_{j},\beta_{j}])\,,\,\mu(\sum_{l}\,[u_{l},\delta_{l}])\,)$.
\end{definition}

Similarly we define fuzzy left ideal, fuzzy right ideal, fuzzy ideal, fuzzy $k$-ideal, fuzzy $h$-ideal of the right operator semiring $R$.
\begin{definition}\cite{Dutta}
Let $S$ be a $\Gamma-$semiring and $L$ be its left operator semiring.\\
For a fuzzy subset $\mu$ of $L$, a fuzzy subset $\mu^{+}$ of $S$ is defined by\\
$\mu^{+}(a)\,=\,\inf_{\alpha\,\in\,\Gamma}\;\;\mu([\,a,\alpha\,])\hspace{1.0cm}where\;a\,\in\,S$.\\
For a fuzzy subset $\sigma$ of $S$, a fuzzy subset $\sigma^{+}\,^{'}$ of $L$ is defined by\\
$\sigma^{+}\,^{'}(\,\sum_{i}\,[x_{i},\alpha_{i}]\,)\,=\,\inf_{a\,\in\,S}\,\sigma(\,\sum_{i}\,x_{i}\alpha_{i}a)\hspace{1.0cm}
where\;\sum_{i}\,[x_{i},\alpha_{i}]\,\in\,L$.
\end{definition}
\begin{definition}\cite{Dutta}
Let $S$ be a $\Gamma-$semiring and $R$ be its right operator semiring.\\
For a fuzzy subset $\mu$ of $R$, a fuzzy subset $\mu^{*}$ of $S$ is defined by\\
$\mu^{*}(a)\,=\,\inf_{\alpha\,\in\,\Gamma}\;\;\mu([\,\alpha,a\,])\hspace{1.0cm}where\;a\,\in\,S$.\\
For a fuzzy subset $\sigma$ of $S$, a fuzzy subset $\sigma^{*}\,^{'}$ of $R$ is defined by\\
$\sigma^{*}\,^{'}(\,\sum_{i}\,[\alpha_{i},x_{i}]\,)\,=\,\inf_{a\,\in\,S}\,\sigma(\,\sum_{i}\,a\alpha_{i}x_{i})\hspace{1.0cm}
where\;\sum_{i}\,[\alpha_{i},x_{i}]\,\in\,R$.
\end{definition}
\begin{definition}\cite{Sardar}
Let $S$ be a $\Gamma-$semiring and $L$ be its left operator semiring and $R$ be its right operator semiring.\\
If there exist an element $\sum_{i}\,[e_{i},\delta_{i}]\,\in\,L$ (respectively $\sum_{j}\,[\gamma_{j},f_{j}]\,\in\,R$)\\s.t.
$\sum_{i}\,e_{i}\delta_{i}a\,=\,a$ (respectively $\sum_{j}\,a\gamma_{j} f_{j}\,=\,a$) for all $a\,\in\,S$ then $S$ is said to have left unity $\sum_{i}\,[e_{i},\delta_{i}]$ (respectively the right unity $\sum_{j}\,[\gamma_{j},f_{j}]$).
\end{definition}
\begin{proposition}\cite{Sardar}
Let $S$ be a $\Gamma$-semiring and $L$ be the left operator semiring of $S$. If  $\sum_{i}\,[e_{i},\delta_{i}]$ is the left unity of $S$, then it is the identity of $L$.
\end{proposition}
\begin{proposition}\cite{Sardar}
Let $S$ be a $\Gamma$-semiring and $R$ be the right operator semiring of $S$. If  $\sum_{j}\,[\gamma_{j},f_{j}]$ is the right unity of $S$, then it is the identity of $R$.
\end{proposition}
Throughout the text, unless otherwise stated explicitly, we consider a $\Gamma$-semiring $S$ which has the left unity, the right unity which implies that the
left operator semiring $L$ and the right operator semiring $R$ of $S$ has identity.
\section{Corresponding Fuzzy Ideals}
\begin{proposition}
Let $S$ be a $\Gamma-$semiring and $L$ be its left operator\\semiring. Then\\
(1) if $\mu$ is a fuzzy left ideal of $L$ then $\mu^{+}$ is a fuzzy left ideal of $S$.\\
(2) if $\mu$ is a fuzzy right ideal of $L$ then $\mu^{+}$ is a fuzzy right ideal of $S$.\\
(3) if $\mu$ is a fuzzy ideal of $L$ then $\mu^{+}$ is a fuzzy ideal of $S$.
\end{proposition}
{\bf Proof.} (1) Let $\mu$ be a fuzzy left ideal of $L$. Let $a,b\in S\,and\,\gamma\in\Gamma$  then\\
$\mu^{+}(a+b)\,=\,\inf_{\alpha\,\in\,\Gamma}\,\mu([a+b,\alpha])$\\
\smallskip\hspace{1.9cm}$\geq\,\inf_{\alpha\,\in\,\Gamma}\,\{\,min\,\{\,\mu([a,\alpha])\,,\,\mu([b,\alpha]\,)\,\}\,\}$\\
\smallskip\hspace{1.9cm}$=\,min\,\{\,\inf_{\alpha\,\in\,\Gamma}\,\mu([a,\alpha])\,,\,\inf_{\alpha\,\in\,\Gamma}\,\mu([b,\alpha])\,\}$\\
\smallskip\hspace{1.9cm}$=\,min\,\{\,\mu^{+}(a)\,,\,\mu^{+}(b)\,\}$\\
$\mu^{+}(a\,\gamma\,b)\;\,=\,\inf_{\alpha\,\in\,\Gamma}\,\mu([a\,\gamma\,b\,,\,\alpha])\hspace{1.5cm}$
where $[a\,\gamma\,b\,,\,\alpha]\,=\,[a\,,\,\gamma].[b\,,\,\alpha]$ \\
\smallskip\hspace{1.9cm}$\geq\,\inf_{\alpha\,\in\,\Gamma}\,\mu([b\,,\,\alpha])\hspace{1cm}$since $\mu$ is a fuzzy left ideal of $L$.\\
\smallskip\hspace{1.9cm}$=\,\mu^{+}(b).$\\
So, $\mu^{+}$ is a fuzzy left ideal of $S$.\\\\
(2) Let $\mu$ be a fuzzy right ideal of $L$. Let $a,b\in S\,and\,\gamma\in\Gamma$.\\
Then by $(1)$ we have $\mu^{+}(a+b)\,\geq\,min\,\{\,\mu^{+}(a)\,,\,\mu^{+}(b)\,\}$\\
Now $\mu^{+}(a\,\gamma\,b)\;\,=\,\inf_{\alpha\,\in\,\Gamma}\,\mu([a\,\gamma\,b\,,\,\alpha])$\\
\smallskip\hspace{2.8cm}$\geq\,\inf_{\alpha\,\in\,\Gamma}\,\mu([a\,,\,\gamma])\hspace{1cm}$since $\mu$ is a fuzzy right ideal of $L$.\\
\smallskip\hspace{2.8cm}$\geq\,\inf_{\alpha\,\in\,\Gamma}\,\mu([a\,,\,\alpha])$\\
\smallskip\hspace{2.8cm}$=\,\mu^{+}(a).$\\
So, $\mu^{+}$ is a fuzzy right ideal of $S$.\\\\
(3) Follows from (1) and (2).
\begin{proposition}
Let $S$ be a $\Gamma-$semiring and $L$ be its left operator\\semiring. Then\\
(1) if$\;\sigma$ is a fuzzy left ideal of $S$ then $\sigma^{+}\,^{'}$ is a fuzzy left ideal of $L$.\\
(2) if$\;\sigma$ is a fuzzy right ideal of $S$ then $\sigma^{+}\,^{'}$ is a fuzzy right ideal of $L$.\\
(3) if$\;\sigma$ is a fuzzy ideal of $S$ then $\sigma^{+}\,^{'}$ is a fuzzy ideal of $L$.
\end{proposition}
{\bf Proof.} (1) Let $\sigma$ be a fuzzy left ideal of $S$. Let $\sum_{i}\,[x_{i},\alpha_{i}]\,,\,\sum_{j}\,[y_{j},\beta_{j}]\,\in\,L$\\
Then $\sigma^{+}\,^{'}(\,\sum_{i}\,[x_{i},\alpha_{i}]\,+\,\sum_{j}\,[y_{j},\beta_{j}]\,)\\
\smallskip\hspace{1cm}=\,\inf_{a\in\,S}\sigma(\,\sum_{i}\,x_{i}\alpha_{i} a\,+\,\sum_{j}\,y_{j}\beta_{j} a\,)\\
\smallskip\hspace{5cm}$since $\sigma$ is a fuzzy left ideal of $S$\\
$\smallskip\hspace{1cm}\geq\,\inf_{a\in\,S}\,\{min\,\{\sigma(\,\sum_{i}\,x_{i}\alpha_{i} a\,)\,,\,\sigma(\sum_{j}\,y_{j}\beta_{j} a\,)\}\}$\\
$\smallskip\hspace{1cm}=\,min\{\,\inf_{a\in\,S}\,\sigma(\,\sum_{i}\,x_{i}\alpha_{i} a\,)\,,\,\inf_{a\in\,S}\,\sigma(\sum_{j}\,y_{j}\beta_{j} a\,)\}\\
\smallskip\hspace{1cm}=\,min\{\,\sigma^{+}\,^{'}(\,\sum_{i}\,[x_{i},\alpha_{i}]\,)\,,\,\sigma^{+}\,^{'}(\,\sum_{j}\,[y_{j},\beta_{j}]\,)\}$.\\
Now   $\,\sigma^{+}\,^{'}(\;\sum_{i}\,[x_{i},\alpha_{i}]\,.\,\sum_{j}\,[y_{j},\beta_{j}]\;)\\
\smallskip\hspace{1cm}=\;\sigma^{+}\,^{'}\,(\,\sum_{i},_{j}[x_{i}\alpha_{i}y_{j}\,,\,\beta_{j}]\,)\\
\smallskip\hspace{1cm}=\,\inf_{a\in S}\,\sigma(\,\sum_{i},_{j}\,(x_{i}\alpha_{i}y_{j})\beta_{j}a)\\
\smallskip\hspace{1cm}=\,\inf_{a\in S}\,\sigma(\,\sum_{i},_{j}\,x_{i}\alpha_{i}(y_{j}\beta_{j}a)\,)\\
\smallskip\hspace{1cm}=\,\inf_{a\in S}\,\sigma(\,\sum_{i}\,x_{i}\alpha_{i}(\sum_{j}\,y_{j}\beta_{j}a)\,)\\
\smallskip\hspace{1cm}\geq\,\inf_{a\in S}\,\{\,min_{i}\,\{\sigma(\,x_{i}\alpha_{i}(\sum_{j}\,y_{j}\beta_{j}a)\,)\}\}\\
\smallskip\hspace{1cm}\geq\,\inf_{a\in S}\,\{\,min_{i}\,\{\sigma(\sum_{j}\,y_{j}\beta_{j}a)\}\}\hspace{1cm}$since $\sigma$ is a fuzzy left ideal of $S$\\
$\smallskip\hspace{1cm}=\,\inf_{a\in S}\,\sigma(\sum_{j}\,y_{j}\beta_{j}a)\\
\smallskip\hspace{1cm}=\,\sigma^{+}\,'(\sum_{j}[y_{j},\beta_{j}])$.\\
So, $\sigma^{+}\,^{'}$ is a fuzzy left ideal of $L$.\\

(2) Suppose $\sigma$ is a fuzzy right ideal of $S$. Let $\sum_{i}\,[x_{i},\alpha_{i}]\,,\,\sum_{j}\,[y_{j},\beta_{j}]\,\in\,L$\\
Then by $(1),\;\sigma^{+}\,^{'}(\,\sum_{i}\,[x_{i},\alpha_{i}]\,+\,\sum_{j}\,[y_{j},\beta_{j}]\,)\,
\geq\,\,min\{\,\sigma^{+}\,^{'}(\,\sum_{i}\,[x_{i},\alpha_{i}]\,)\,,\,\sigma^{+}\,^{'}(\,\sum_{j}\,[y_{j},\beta_{j}]\,)\}$.\\
Now  $\sigma^{+}\,^{'}(\;\sum_{i}\,[x_{i},\alpha_{i}]\,.\,\sum_{j}\,[y_{j},\beta_{j}]\;)\\
\smallskip\hspace{1cm}=\;\sigma^{+}\,^{'}\,(\,\sum_{i},_{j}[x_{i}\alpha_{i}y_{j}\,,\,\beta_{j}]\,)\\
\smallskip\hspace{1cm}=\,\inf_{a\in S}\,\sigma(\,\sum_{i},_{j}\,(x_{i}\alpha_{i}y_{j})\beta_{j}a\,)\\
\smallskip\hspace{1cm}=\,\inf_{a\in S}\,\sigma(\,\sum_{j}(\sum_{i}\,x_{i}\alpha_{i}y_{j})\beta_{j}a\,)\\
\smallskip\hspace{1cm}\geq\,\inf_{a\in S}\,\{min_{j}\,\{\sigma(\sum_{i}\,x_{i}\alpha_{i}y_{j})\beta_{j}a\,)\}\}\\
\smallskip\hspace{1cm}\geq\,\inf_{a\in S}\,\{min_{j}\,\{\sigma(\sum_{i}\,x_{i}\alpha_{i}y_{j}\,)\}\}$\hspace{0.5cm} since $\sigma$ is a fuzzy right ideal of $S$\\
$\smallskip\hspace{1.1cm}=\,min_{j}\,\{\sigma(\sum_{i}\,x_{i}\alpha_{i}y_{j}\,)\}\\
\smallskip\hspace{1cm}=\,\sigma^{+}\,^{'}\,(\,\sum_{i}\,[x_{i},\alpha_{i}]\,)$.\\
So, $\sigma^{+}\,^{'}$ is a fuzzy right ideal of $L$.\\\\
(3) Follows from (1) and (2).\\
\begin{Theorem}
Let $S$ be a $\Gamma$-semiring and $L$ be its left operator semiring.\\
Then there exist an inclusion preserving bijection $\sigma\,\rightarrow\,\sigma^{+}\,^{'}$ between the lattices of all fuzzy right ideals (respectively fuzzy ideals)
of $S$ and the lattices of all fuzzy right ideals (respectively fuzzy ideals) of $L$.\\
Where  $\sigma$ denotes a fuzzy right ideals (respectively fuzzy ideals) of $S$.
\end{Theorem}
{\bf Proof.} We first consider the case for lattices of all fuzzy right ideals.\\
Let $\sigma$ be a fuzzy right ideal of $S$. Then for all $a\,\in S$,\\
$(\sigma^{+}\,^{'})^{+}\,(a)\,=\,\inf_{\alpha\in\Gamma}\,\sigma^{+}\,^{'}(\,[a,\alpha]\,)=\,\inf_{\alpha\,\in\Gamma}\,\inf_{b\,\in S}\,\sigma(a\alpha b)\;\;$..................(1)\\
Since $S$ has the right unity say, $\sum_{i}\,[\gamma_{i},f_{i}]$\\
So, for all $a\,\in\,S,\;\;a\,=\,\sum_{i}\,a\gamma_{i}f_{i}$\\
Therefore, $\sigma(a)\,=\,\sigma(\sum_{i}\,a\gamma_{i}f_{i})\,\geq\,min_{i}\,\sigma(a\gamma_{i}f_{i})\,\geq\,\inf_{\alpha\,\in\,\Gamma}\inf_{b\,\in\, S}\,\sigma(a\alpha b)$\\
So, by (1) we have $(\sigma^{+}\,^{'})^{+}\,(a)\,\leq\,\sigma(a)\;\;$............................................(2)\\
\smallskip\hspace{.5cm}Again, since $\sigma$ is a fuzzy right ideal of $S$, therefore $\sigma(a\alpha b)\,\geq\,\sigma(a)$ for all $a,b\,\in S,\,\alpha\,\in\Gamma$\\
So, by (1) we have $(\sigma^{+}\,^{'})^{+}\,(a)\,\geq\,\inf_{\alpha\,\in\Gamma}\,\inf_{b\,\in S}\,\sigma(a)\,=\,\sigma(a)\;\;$..........(3)\\
By (2) and (3) we have $(\sigma^{+}\,^{'})^{+}\,(a)\,=\,\sigma(a)$. Hence the mapping $\sigma\,\rightarrow\,\sigma^{+}\,^{'}$ is one-to-one.\\

Now let $\mu$ be a fuzzy right ideal of $L$. Then for all $\sum_{i}\,[x_{i},\alpha_{i}]\,\in\,L$,\\
$(\mu^{+})^{+}\,^{'}\,(\sum_{i}\,[x_{i},\alpha_{i}])\,=\,\inf_{a\in S}\,\mu^{+}\,(\sum_{i}\,x_{i}\alpha_{i}a)\,=\,\inf_{a\in S}\,\inf_{\alpha\in\Gamma}\,\mu(\,[\,\sum_{i}x_{i}\alpha_{i}a\,,\alpha\,]\,)\\
=\,\inf_{a\in S}\,\inf_{\alpha\in\Gamma}\,\mu(\,\sum_{i}\,[x_{i}\alpha_{i}a,\alpha]\,)$\hspace{1.2cm} since $[x+y,\alpha]\,=\,[x,\alpha]+[y,\alpha]$\\
$=\,\inf_{a\in S}\,\inf_{\alpha\in\Gamma}\,\mu(\,(\sum_{i}\,[x_{i},\alpha_{i}]).[a,\alpha]\,)\;\;$......................................(4)\\
Since $L$ has the identity say, $\sum_{j}\,[e_{j},\delta_{j}]$\\
So, $\sum_{i}\,[x_{i},\alpha_{i}]\,.\sum_{j}\,[e_{j},\delta_{j}]\,=\,\sum_{i}\,[x_{i},\alpha_{i}]$\\
Therefore, $\mu(\,\sum_{i}\,[x_{i},\alpha_{i}]\,)\,=\,\mu(\,\sum_{i}\,[x_{i},\alpha_{i}]\,.\sum_{j}\,[e_{j},\delta_{j}]\,)\,=\,\mu(\,\sum_{i},_{j}\,[x_{i}\alpha_{i}e_{j},\delta_{j}]\,)\\
=\,\mu(\,\sum_{j}\,[\sum_{i}x_{i}\alpha_{i}e_{j},\delta_{j}]\,)\,\geq\,min_{j}\,\mu(\sum_{i}\,[x_{i}\alpha_{i}e_{j},\delta_{j}])\\
\smallskip\hspace{7cm}$ since $\mu$ is a fuzzy right ideal of $L$\\
$\geq\,\inf_{a\in S}\inf_{\alpha\in\Gamma}\,\mu(\,\sum_{i}\,[x_{i}\alpha_{i}a,\alpha]\,)$\\
So, by (4) we have $(\mu^{+})^{+}\,^{'}\,(\sum_{i}\,[x_{i},\alpha_{i}])\,\leq\,\mu(\,\sum_{i}\,[x_{i},\alpha_{i}]\,)\;\;$..........(5)\\
\smallskip\hspace{.5cm}Again, since $\mu$ is a fuzzy right ideal of $L$, therefore $\mu(\,(\sum_{i}\,[x_{i},\alpha_{i}]).[a,\alpha]\,)\,\geq\,\mu(\,\sum_{i}\,[x_{i},\alpha_{i}]\,)$ where $\sum_{i}\,[x_{i},\alpha_{i}]\,,[a,\alpha]\,\in\,L$\\
So, by (4) we have $(\mu^{+})^{+}\,^{'}\,(\sum_{i}\,[x_{i},\alpha_{i}])\,\geq\,\inf_{a\in S}\inf_{\alpha\in\Gamma}\,\mu(\,\sum_{i}\,[x_{i},\alpha_{i}]\,)\\
\smallskip\hspace{6.9cm}=\,\mu(\,\sum_{i}\,[x_{i},\alpha_{i}]\,)\;$.........(6)\\
By (5) and (6) we have $(\mu^{+})^{+}\,^{'}\,(\sum_{i}\,[x_{i},\alpha_{i}])\,=\,\mu(\,\sum_{i}\,[x_{i},\alpha_{i}]\,)$.\\
Hence the mapping $\sigma\,\rightarrow\,\sigma^{+}\,^{'}$ is onto. Thus the mapping is bijective.\\

Now let $\sigma$ and $\mu$ be two fuzzy right ideals of $S$ s.t. $\sigma\leq\mu$. Then for all $\sum_{i}\,[x_{i},\alpha_{i}]\,\in\,L$.
\smallskip\hspace{.5cm}$\sigma^{+}\,^{'}\,(\,\sum_{i}\,[x_{i},\alpha_{i}]\,)\,=\,\inf_{a\in S}\,\sigma(\,\sum_{i}\,x_{i}\alpha_{i}a\,)\,\leq\,\inf_{a\in S}\,\mu(\,\sum_{i}\,x_{i}\alpha_{i}a\,)\\=\,\mu^{+}\,^{'}(\,\sum_{i}\,[x_{i},\alpha_{i}]\,)$. Hence the said mapping is order preserving.\\
Also let $\eta\,,\,\nu$ be two fuzzy right ideals of $L$ s.t. $\eta\leq\nu$. Then for all $a\,\in\,S$\\
$\eta^{+}(a)\,=\,\inf_{\alpha\in\Gamma}\,\eta(\,[a,\alpha]\,)\,\leq\,\inf_{\alpha\in\Gamma}\,\nu(\,[a,\alpha]\,)\,=\,\nu^{+}(a)$.\\
So, inverse of the said mapping is also order preserving.\\
Hence for the lattices of all fuzzy right ideals the theorem is proved.\\
Similar proof for the case for the lattices of all fuzzy two-sided ideals or simply fuzzy ideals. Hence the theorem is proved.\\

Now we obtain analogous results of Propositions 3.1, 3.2 and Theorem 3.3 for the right operator semiring $R$ of a $\Gamma$-semiring $S$. We only state the results as the proof is similar as in case of the left operator semiring $L$.\\
\begin{proposition}
Let $S$ be a $\Gamma-$semiring and $R$ be its right operator\\semiring. Then\\
(1) if $\mu$ is a fuzzy left (respectively right) ideal of $R$ then $\mu^{*}$ is a fuzzy left (respectively right) ideal of $S$.\\
(2) if $\mu$ is a fuzzy ideal of $R$ then $\mu^{*}$ is a fuzzy ideal of $S$.
\end{proposition}
\begin{proposition}
Let $S$ be a $\Gamma-$semiring and $R$ be its right operator\\semiring. Then\\
(1) if $\sigma$ is a fuzzy left (respectively right) ideal of $S$ then $\sigma^{*}\,^{'}$ is a fuzzy left (respectively right) ideal of $R$.\\
(2) if $\sigma$ is a fuzzy ideal of $S$ then $\sigma^{*}\,^{'}$ is a fuzzy ideal of $R$.
\end{proposition}
\begin{Theorem}
Let $S$ be a $\Gamma-$semiring and $R$ be its right operator semiring.
Then there exist an inclusion preserving bijection $\sigma\,\rightarrow\,\sigma^{*}\,^{'}$ between the lattices of all fuzzy left ideals (respectively fuzzy ideals)
of $S$ and the lattices of all fuzzy left ideals (respectively fuzzy ideals) of $R$.\\
Where  $\sigma$ denotes a fuzzy left ideal (respectively fuzzy ideal) of $S$.
\end{Theorem}
\begin{proposition}
Let $S$ be a $\Gamma-$semiring and $L$ be its left operator semiring.\\
Then (1) if $\sigma$ is a fuzzy $k-$ideal of $S$ then $\sigma^{+}\,^{'}$ is a fuzzy $k-$ideal of $L$.\\
\smallskip\hspace{1cm}(2) if $\mu$ is a fuzzy $k-$ideal of $L$ then $\mu^{+}$ is a fuzzy $k-$ideal of $S$.
\end{proposition}
{\bf Proof.} (1) Let $\sigma$ be a fuzzy $k-$ideal of $S$. Then by proposition 3.2\\
$\sigma^{+}\,^{'}$ is a fuzzy ideal of $L$.  Now, let $\sum_{i}\,[x_{i},\alpha_{i}]\,,\,\sum_{j}\,[y_{j},\beta_{j}]\,\in\,L$,\\
then $\sum_{i}\,x_{i}\alpha_{i}a\,,\,\sum_{j}\,y_{j}\beta_{j}a\,\in\,S$ for all $a\in\,S$.\\
As $\sigma$ is a fuzzy $k-$ideal of $S$ then\\ $\sigma(\sum_{i}\,x_{i}\alpha_{i}a)\,\geq\,min(\,\sigma(\sum_{i}\,x_{i}\alpha_{i}a+\sum_{j}\,y_{j}\beta_{j}a)\,,\,\sigma(\sum_{j}\,y_{j}\beta_{j}a)\,)$.\\
So, $\inf_{a\in\,S}\,\sigma(\sum_{i}\,x_{i}\alpha_{i}a)\,\geq\,min(\,\inf_{a\in\,S}\,\sigma(\sum_{i}\,x_{i}\alpha_{i}a+\sum_{j}\,y_{j}\beta_{j}a)\,,\,\inf_{a\in\,S}\,\sigma(\sum_{j}\,y_{j}\beta_{j}a)\,)$.\\
This implies that $\sigma^{+}\,^{'}(\sum_{i}\,[x_{i},\alpha_{i}])\,\geq\,min(\,\sigma^{+}\,^{'}(\sum_{i}\,[x_{i},\alpha_{i}]+\sum_{j}\,[y_{j},\beta_{j}])\,,\,\sigma^{+}\,^{'}(\sum_{j}\,[y_{j}\beta_{j}])\,)$.\\
Hence $\sigma^{+}\,^{'}$ is a fuzzy $k-$ideal of $L$.\\

(2) Let $\mu$ be a fuzzy $k-$ideal of $L$. Then by proposition 3.1\\$\mu^{+}$ is a fuzzy ideal of $S$.\\
Now, let $x,y\in\,S$,  then $[x,\alpha]\,,\,[y,\alpha]\in\,L$ for all $\alpha\in\,\Gamma$.\\
As $\mu$ is a fuzzy $k-$ideal of $L$, then\\
$\mu([x,\alpha])\,\geq\,min(\,\mu([x,\alpha]+[y,\alpha])\,,\,\mu([y,\alpha])\,)\,=\,min(\,\mu([x+y,\alpha])\,,\,\mu([y,\alpha])\,)$\\
So, $\inf_{\alpha\in\,\Gamma}\,\mu([x,\alpha])\,\geq\,min(\,\inf_{\alpha\in\,\Gamma}\,\mu([x+y,\alpha])\,,\,\inf_{\alpha\in\,\Gamma}\,\mu([y,\alpha])\,)$.\\
This implies that $\mu^{+}(x)\,\geq\,min(\,\mu^{+}(x+y)\,,\,\mu^{+}(y)\,)$.\\
Hence $\mu^{+}$ is a fuzzy $k-$ideal of $S$.
\begin{proposition}
Let $S$ be a $\Gamma-$semiring and $L$ be its left operator semiring.\\
Then (1) if $\sigma$ is a fuzzy $h-$ideal of $S$ then $\sigma^{+}\,^{'}$ is a fuzzy $h-$ideal of $L$.\\
\smallskip\hspace{1cm}(2) if $\mu$ is a fuzzy $h-$ideal of $L$ then $\mu^{+}$ is a fuzzy $h-$ideal of $S$.
\end{proposition}
{\bf Proof.} (1) Let $\sigma$ be a fuzzy $h-$ideal of $S$. Then by proposition 3.2\\$\sigma^{+}\,^{'}$ is a fuzzy ideal of $L$.\\
Now, let $\sum_{i}\,[x_{i},\alpha_{i}],\;\sum_{j}\,[y_{j},\beta_{j}],\;\sum_{k}\,[z_{k},\gamma_{k}],\;\sum_{l}\,[u_{l},\eta_{l}]\,\in L$ such that\\
$\sum_{i}\,[x_{i},\alpha_{i}]+\sum_{j}\,[y_{j},\beta_{j}]+\sum_{k}\,[z_{k},\gamma_{k}]\,=\,\sum_{l}\,[u_{l},\eta_{l}]+\sum_{k}\,[z_{k},\gamma_{k}]$.\\
Then $\sum_{i}\,x_{i}\alpha_{i}a,\,\sum_{j}\,y_{j}\beta_{j}a,\,\sum_{k}\,z_{k}\gamma_{k}a,\,\sum_{l}\,u_{l}\eta_{l}a\,\in S$ for all $a\in S$ and\\ $\sum_{i}\,x_{i}\alpha_{i}a+\sum_{j}\,y_{j}\beta_{j}a+\sum_{k}\,z_{k}\gamma_{k}a\,=\,\sum_{l}\,u_{l}\eta_{l}a+\sum_{k}\,z_{k}\gamma_{k}a$.\\
As $\sigma$ is a fuzzy $h-$ideal of $S$ then\\ $\sigma(\sum_{i}\,x_{i}\alpha_{i}a)\,\geq\,min(\,\sigma(\sum_{j}\,y_{j}\beta_{j}a)\,,\,\sigma(\sum_{l}\,u_{l}\eta_{l}a)\,)$ for all $a\in S$.\\
So, $\inf_{a\in\,S}\,\sigma(\sum_{i}\,x_{i}\alpha_{i}a)\,\geq\,min(\,\inf_{a\in\,S}\,\sigma(\sum_{j}\,y_{j}\beta_{j}a)\,,\,\inf_{a\in
\,S}\,\sigma(\sum_{l}\,u_{l}\eta_{l}a)\,)$.\\
This implies that $\sigma^{+}\,^{'}(\sum_{i}\,[x_{i},\alpha_{i}])\,\geq\,min(\,\sigma^{+}\,^{'}(\sum_{j}\,[y_{j},\beta_{j}])\,,\,\sigma^{+}\,^{'}(\sum_{l}\,[u_{l},\eta_{l}])\,)$.\\
Hence $\sigma^{+}\,^{'}$ is a fuzzy $h-$ideal of $S$.\\

(2) Let $\mu$ be a fuzzy $h-$ideal of $L$. Then by proposition 3.1 $\mu^{+}$ is a fuzzy ideal of $S$.\\
Now, let $x,z,y_{1},y_{2}\in S$ such that $x+y_{1}+z\,=\,y_{2}+z$.\\
Then for all $\alpha\in \Gamma$, we have $[x,\alpha]\,,\,[z,\alpha]\,,\,[y_{1},\alpha]\,,\,[y_{2},\alpha]\in L$ and\\
$[x+y_{1}+z,\alpha]\,=\,[y_{2}+z,\alpha]$ i.e. $[x,\alpha]+[y_{1},\alpha]+[z,\alpha]\,=\,[y_{2},\alpha]+[z,\alpha]$.\\
As $\mu$ is a fuzzy $h-$ideal of $L$ then $\mu([x,\alpha])\,\geq\,min(\,\mu([y_{1},\alpha]),\mu([y_{2},\alpha])\,)$ for all $\alpha\in \Gamma$.\\
So, $\inf_{\alpha\in\,\Gamma}\,\mu([x,\alpha])\,\geq\,min(\,\inf_{\alpha\in\,\Gamma}\,\mu([y_{1},\alpha])\,,\,\inf_{\alpha\in\,\Gamma}\,\mu([y_{2},\alpha])\,)$.\\
This implies that $\mu^{+}(x)\,\geq\,min(\,\mu^{+}(y_{1})\,,\,\mu^{+}(y_{2})\,)$.\\
Hence $\mu^{+}$ is a fuzzy $h-$ideal of $S$.\\

Now we obtain analogous results of Propositions 3.7, 3.8 for the right operator semiring $R$ of a $\Gamma$-semiring $S$. We only state the results as the proof is similar as in case of the left operator semiring $L$.
\begin{proposition}
Let $S$ be a $\Gamma-$semiring and $R$ be its right operator semiring.\\
Then (1) if $\sigma$ is a fuzzy $k-$ideal of $S$ then $\sigma^{*}\,^{'}$ is a fuzzy $k-$ideal of $R$.\\
\smallskip\hspace{1cm}(2) if $\mu$ is a fuzzy $k-$ideal of $R$ then $\mu^{*}$ is a fuzzy $k-$ideal of $S$.
\end{proposition}
\begin{proposition}
Let $S$ be a $\Gamma-$semiring and $R$ be its left operator semiring. Then\\
(1) if $\sigma$ is a fuzzy $h-$ideal of $S$ then $\sigma^{*}\,^{'}$ is a fuzzy $h-$ideal of $R$.\\
(2) if $\mu$ is a fuzzy $h-$ideal of $R$ then $\mu^{*}$ is a fuzzy $h-$ideal of $S$.
\end{proposition}
\begin{remark}
Propositions 3.7,3.8,3.9 and 3.10 are also valid for fuzzy one-sided $k-$ideals and $h-$ideals.
\end{remark}
We now only state the following theorems as the proof is analogous to that of the Theorem 3.3.
\begin{Theorem}
Let $S$ be a $\Gamma$-semiring and $L$ be its left operator semiring.\\
Then there exist an inclusion preserving bijection $\sigma\,\rightarrow\,\sigma^{+}\,^{'}$ between the set of all fuzzy $k-$ideals (respectively fuzzy $h-$ideals)
of $S$ and the set of all fuzzy $k-$ideals (respectively fuzzy $h-$ideals) of $L$.\\
Where  $\sigma$ denotes a fuzzy $k-$ideal (respectively fuzzy $h-$ideal) of $S$.
\end{Theorem}
\begin{Theorem}
Let $S$ be a $\Gamma$-semiring and R be its right operator semiring.\\
Then there exist an inclusion preserving bijection $\sigma\,\rightarrow\,\sigma^{*}\,^{'}$ between the set of all fuzzy $k-$ideals (respectively fuzzy $h-$ideals)
of $S$ and the set of all fuzzy $k-$ideals (respectively fuzzy $h-$ideals) of $R$.\\
Where  $\sigma$ denotes a fuzzy $k-$ideal (respectively fuzzy $h-$ideal) of $S$.
\end{Theorem}
\begin{remark}
Theorem 3.12 is also valid for fuzzy right $k-$ideals and fuzzy right $h-$ideals and Theorem 3.13 is also valid for fuzzy left $k-$ideals and fuzzy left $h-$ideals.
\end{remark}

\end{document}